\documentclass[12pt,twoside]{amsart}
\usepackage{systeme}
\usepackage{graphicx}
\usepackage{tikz}
\usepackage{tikz-cd}
\usetikzlibrary{decorations.markings}
\usetikzlibrary{arrows,shapes,positioning}
\usetikzlibrary{topaths}
\tikzstyle{vertex}=[circle, draw, inner sep=0pt, minimum size=4pt]
\usetikzlibrary{positioning}

\usepackage{amssymb,amscd}
\usepackage{mathrsfs}
\usepackage{array,float}
\usepackage[all]{xy}
\usepackage{enumerate}

\theoremstyle{plain}
\newtheorem{thm}{Theorem}[section]
\newtheorem{thmx}{Theorem}

\newtheorem{lem}[thm]{Lemma}
\newtheorem{pro}[thm]{Proposition}
\newtheorem{cor}[thm]{Corollary}

\theoremstyle{remark}
\newtheorem{rem}[thm]{Remark}

\theoremstyle{definition} 
\newtheorem{definition}[thm]{Definition}

\numberwithin{equation}{section}

\renewcommand{\phi}{\varphi}

\begin{document}

\author{Ilir Snopce}
\address{Universidade Federal do Rio de Janeiro\\
 Instituto de Matem\'atica \\
21941-909 Rio de Janeiro, RJ \\ Brasil }
\email{ilir@im.ufrj.br}

\title[Frattini-resistance and direct products] {Frattini-resistant  direct products of pro-$p$ groups}

\author{Slobodan Tanushevski} 
\address{Universidade Federal Fluminense\\
Instituto de Matem\'atica e Estat\'istica\\
24210-201  Nit\'eroi, RJ \\ Brasil }
\email{stanushevski@id.uff.br}

\subjclass[2010]{20E18, 12F10, 22E20}

\thanks{The first  author is partially supported by FAPERJ and CNPq.}

\begin{abstract} 
A pro-$p$ group $G$ is called strongly Frattini-resistant if the function $H \mapsto \Phi(H)$, from the poset of all closed subgroups of $G$ into itself, is a poset embedding.  
Frattini-resistant pro-$p$ groups appear naturally in Galois theory. Indeed, every maximal pro-$p$ Galois group over a field that contains a primitive $p$th root of unity (and also contains $\sqrt{-1}$ if $p=2$) is strongly Frattini-resistant.
Let $G_1$ and $G_2$ be non-trivial pro-$p$ groups. We prove that $G_1 \times G_2$ is strongly Frattini-resistant if and only if one of the direct factors $G_1$ or $G_2$ is torsion-free abelian and the other one has the property that all of its closed 
subgroups have torsion-free abelianization. As a corollary we obtain a group theoretic proof of a result of  Koenigsmann on maximal pro-$p$ Galois groups that admit a non-trivial decomposition as a direct product. 
In addition, we give an example of a group that is not strongly Frattini-resistant, but has the property that 
its Frattini-function defines an order self-embedding of the poset of all topologically finitely generated subgroups.      
\end{abstract}

\maketitle
\section{Introduction}

Given a pro-$p$ group $G$, let $\mathcal{S}(G)$ denote the poset of closed subgroups of $G$, and let $\mathcal{S}_{f.g.}(G)$ be the subposet consisting of all topologically finitely generated subgroups. The Frattini function ${\Phi: \mathcal{S}(G) \to \mathcal{S}(G)}$, $\Phi(H)=\overline{H^p[H, H]}$, is an order-preserving self-mapping of $\mathcal{S}(G)$. Since the Frattini subgroup of a (topologically) finitely generated pro-$p$ group is also finitely generated, it follows that the Frattini function restricts to a  self-mapping of $\mathcal{S}_{f.g.}(G)$.

In \cite{SnoTan}, we introduced the following concept:

\begin{definition}
A pro-$p$ group $G$ is called Frattini-resistant if the (restricted) Frattini function $\Phi:\mathcal{S}_{f.g.}(G) \to \mathcal{S}_{f.g.}(G)$ is an order embedding, i.e., for all ${H, K \in \mathcal{S}_{f.g.}(G)}$, 
\[\Phi(H)\leq \Phi(K) \, \Leftrightarrow \, H \leq K.\] 
If the stronger condition that $\Phi: \mathcal{S}(G) \to \mathcal{S}(G)$ is an order embedding is satisfied, then we say that $G$ is strongly Frattini-resistant. 
\end{definition}

Following \cite{Wurfel}, we call a pro-$p$ group $G$ \emph{absolutely torsion-free} if every subgroup of $G$ has torsion-free (and thus free abelian) abelianization. 
Every absolutely torsion-free pro-$p$ group is strongly Frattini-resistant (cf. \cite[Proposition~4.10]{SnoTan}). 

The main result of this paper provides a necessary and sufficient condition for a direct product of  non-trivial pro-$p$ groups to be (strongly) Frattini-resistant. (It can be seen as a generalization of \cite[Proposition~4]{Wurfel}.)
 
\begin{thmx}
\label{main theorem}
Let $p$ be any prime, and let $G_1$ and $G_2$ be non-trivial pro-$p$ groups.  The following assertions hold:
\begin{enumerate}[(i)]
\item $G_1 \times G_2$ is a Frattini-resistant pro-$p$ group if and only if  one of the groups $G_1$ or $G_2$ is torsion-free abelian and 
the other one is locally absolutely torsion-free, in which case, $G_1 \times G_2$ is locally absolutely torsion-free.
\item $G_1 \times G_2$ is a strongly Frattini-resistant pro-$p$ group if and only if  one of the groups $G_1$ or $G_2$ is torsion-free abelian and 
the other one is absolutely torsion-free, in which case, $G_1 \times G_2$ is absolutely torsion-free.
\end{enumerate}
\end{thmx}



Let $k$ be a field containing a primitive $p$th root of unity; if $p=2$ assume in addition that $\sqrt{-1} \in k$.
Let $k(p)$ denote the maximal $p$-extension of $k$, that is,  the compositum of all finite Galois $p$-extensions of $k$ inside a fixed seperable closure of $k$.  
The Galois group $G_k(p):=\mathrm{Gal}(k(p)/k)$ is called the maximal pro-$p$ Galois group of $k$. By \cite[Theorem 7.1 and Theorem 7.2]{SnoTan}, $G_k(p)$  is a strongly Frattini-resistant pro-$p$ group.
In fact, we do not know of any example of a strongly Frattini-resistant pro-$p$ group that can not be realized as a maximal pro-$p$ Galois group over a field satisfying the above mentioned condition. Moreover, numerous results on maximal pro-$p$ Galois groups 
can be derived as rather elementary group-theoretic consequences of Frattini-resistance, as we show in  \cite{SnoTan}.

As a corollary of Theorem~\ref{main theorem}, we obtain yet another substantial result on maximal pro-$p$ Galois groups.

\begin{cor}[{\cite[Proposition~2.2]{Koen}}
\label{product of Galois groups}]
Let $k$ be a field containing a primitive $p$th root of unity, and let $\mu_{\infty}$ denote the set of all roots of unity of order a power $p$ inside $k(p)$.
If $G_k(p) \cong G_1 \times G_2$ with $G_1$ and $G_2$ non-trivial pro-$p$ groups, then $\mu_{\infty} \subseteq k$ and 
one of the direct factors $G_1$ or $G_2$ is torsion-free abelian and the other one is absolutely torsion-free; in particular, $G_k(p)$ is absolutely torsion-free.
\end{cor}

Theorem~\ref{main theorem} implies a similar result for $1$-smooth pro-$p$ groups. (For the definition and results on $1$-smooth pro-$p$ groups we refer the reader to \cite{CleFlo17}, \cite{QuWe18} and \cite[Section~ 7]{SnoTan}.)

\begin{cor}
\label{smooth}
Let $G_1$ and $G_2$ be non-trivial pro-$p$ groups.  
There exists a homomorphism $\theta: G_1 \times G_2 \to 1 + p\mathbb{Z}_p$  (with $\textrm{Im}(\theta) \leq 1 + 4\mathbb{Z}_2$ if $p=2$)   
such that $(G_1 \times G_2, \theta)$ is a $1$-smooth cyclotomic pro-$p$ pair  if and only if one of the groups $G_1$ or $G_2$ is torsion-free abelian and 
the other one is absolutely torsion-free, in which case $\theta$ must be the trivial homomorphism. 
\end{cor}

A finitely generated Frattini-resistant pro-$p$ group is necessarily strongly Frattini-resistant (\cite[Corollary 4.5]{SnoTan}).  
Theorem~\ref{main theorem} provides means for constructing tangible examples of infinitely generated Frattini-resistant pro-$p$ groups that are not strongly Frattini-resistant. 

\begin{cor}
\label{not strongly Frattini-resistant}
Let $n$ be a positive integer, and let $G$ be the Demushkin pro-$p$ group of infinite rank defined by the presentation
\[G:=\langle x_1, x_2, x_3 \ldots \mid x_1^{p^n}[x_1, x_2][x_3, x_4][x_5, x_6]\cdots \rangle.\] 
The pro-$p$ group $G \times \mathbb{Z}_p$ is Frattini-resistant but not strongly Frattini-resistant.
\end{cor}    

\textbf{Notation:}
When discussing pro-$p$ groups, all group theoretic notions are taken in the appropriate sense for topological groups; thus subgroups are assumed to be closed, all homomorphisms are
continuous, and  generators are always understood to be topological generators. Let $G$ be a pro-$p$ group and $x, y \in G$. We use the following
notation: $d(G)$ is the smallest cardinality of a generating set for $G$; $[x, y]=x^{-1}y^{-1}xy$; the $n$th term of the lower
central series of $G$ is denoted by $\gamma_{n}(G)$, with the exception
of the commutator subgroup, which is always denoted by $[G, G]$; we write $G^{ab}$ for the abelianization of $G$. 
\section{Preliminaries}

A more manageable  rendition of Frattini-resistance is provided by the following  

\begin{pro}[{\cite[Proposition~4.2]{SnoTan}}]
Let $G$ be a pro-$p$ group. The following assertions hold: 
\begin{enumerate}[(i)]
\item $G$ is Frattini-resistant if and only if for all $H \in \mathcal{S}_{f.g.}(G)$ and $x \in G$,
\[x^p \in \Phi(H) \Rightarrow x \in H.\]
\item $G$ is strongly Frattini-resistant if and only if for all $H \in \mathcal{S}(G)$ and $x \in G$,
\[x^p \in \Phi(H) \Rightarrow x \in H.\]
\end{enumerate}
\end{pro}  

\medskip

In \cite{SnoTan}, we defined a pro-$p$ group $G$ to be \emph{commutator-resistant} if for every 
${H \in \mathcal{S}_{f.g.}(G)}$ and every $x\in H$, $x^p \in [H, H] \Rightarrow x \in \Phi(H)$. 
If this condition holds for all subgroups $H$ of $G$, then we say that $G$ is \emph{strongly commutator-resistant}.

Perhaps a more instructive way of thinking about commutator-resistance is afforded by the following equivalence:
A pro-$p$ group $G$ is commutator-resistant if and only if  for every $H \in \mathcal{S}_{f.g.}(G)$,
$H^{ab}$ does not contain a cyclic direct factor of order $p$ (cf. {\cite[Corollary~4.12]{SnoTan}}). 

By {\cite[Proposition~4.10]{SnoTan}}, every (strongly) commutator-resistant pro-$p$ group is (strongly) Frattini-resistant. However, for every odd prime $p$, there exist finitely generated strongly Frattini-resistant pro-$p$ groups that are not commutator-resistant. (For explicit examples in the class of Demushkin groups see {\cite[Theorem~6.3]{SnoTan}}.) It turns out that for pro-$2$ groups the two concepts coincide.

\begin{thm}
\label{commutator-resistant}
Let $G$ be a pro-$2$ group. The following assertions hold:
\begin{enumerate}[(i)]
\item $G$ is Frattini-resistant if and only if it is commutator-resistant.
\item $G$ is strongly Frattini-resistant if and only if it is strongly commutator-resistant.  
\end{enumerate}
\end{thm}
\begin{proof}
$(i)$ Suppose to the contrary that there is a Frattini-resistant pro-$2$ group $G$ that is not commutator-resistant. Thus there exist a finitely generated subgroup $H$ of $G$ and an element $h \in H \setminus \Phi(H)$ such that $h^2 \in [H, H]$.
Let $M$ be a maximal subgroup of $H$ that does not contain $h$. Then there is $m \in M$ such that 
\[h^2 \Phi(M)=[h, m]\Phi(M).\] 
(Indeed, for every $w \in [H, H]$, we have $w\Phi(M)=[h, x]\Phi(M)$ for some $x \in M $.)  
It follows that 
\[(hm)^2\Phi(M)=h^2m^2[m, h][[m, h], m]\Phi(M)=h^2[m, h]\Phi(M)=\Phi(M),\]
that is, $(hm)^2 \in \Phi(M)$. Since $hm \notin M$, we obtain a contradiction with the Frattini-resistance of $G$.

 Since in the proof of $(i)$ we did not use the fact that $H$ is finitely generated, the same argument yields $(ii)$.



\end{proof}

A finitely generated pro-$p$ group $G$ is absolutely torsion-free if  $H^{ab}$ is torsion-free for every 
$H \in \mathcal{S}_{f.g.}(G)$. Indeed, let $G$ be a finitely generated pro-$p$ group having the property that
all of its finitely generated subgroups have torsion free abelianization. Thus, in particular, 
every open subgroup of $G$ has torsion free abelianization. Now let $H$ be any subgroup of $G$, and  
assume that $H^{ab}$ has a non-trivial element of finite order. Then there is $h \in H \setminus \Phi(H)$ with 
$h[H, H]$ an element of finite order in $H^{ab}$. Choose an open subgroup $U$ of $G$ such that 
$H \leq U$ and $h \notin [U, U]$ (cf. {\cite[Proposition~2.8.9]{Ribes1}}). Since $[H, H] \leq [U, U]$, it follows that $h[U, U]$ is a non-trivial element of finite order in $U^{ab}$, which contradicts our assumption.
     
It follows from the observation made in the previous paragraph that the designation `locally absolutely torsion-free' is nothing but an abbreviation for  `all finitely generated subgroups have torsion-free abelianization'.

\begin{pro}
\label{Demushkin of inf rank}
Let $n$ be a positive integer. The Demushkin pro-$p$ group of infinite rank defined by the presentation
\[G:=\langle x_1, x_2, x_3 \ldots \mid x_1^{p^n}[x_1, x_2][x_3, x_4][x_5, x_6]\cdots \rangle\]
is locally absolutely torsion-free but not absolutely torsion-free. 
\end{pro}
\begin{proof}
Since $x_1^{p^n} \in [G, G]$, it follows that $G^{ab}$ has a non-trivial element of finite order, and therefore, $G$ is not absolutely torsion-free.

\medskip

For each positive integer $i$, let $G_i$ be the Demushkin group (of finite rank) defined by the presentation
\[G_i:= \langle x_1, x_2, \ldots, x_{2i} \mid x_1^{p^n}[x_1, x_2][x_3, x_4] \cdots [x_{2i-1}, x_{2i}] \rangle.\] 
For all $i, j \in \mathbb{N}$ with $i \leq j$, there is a homomorphism $\phi_{i, j}:G_j \to G_i$ defined by $\phi_{i, j}(x_k)=x_k$ for $1 \leq k \leq 2i$ and $\phi_{i, j}(x_k)=1$ for $k>2i$.
In the same vein, for each $i \in \mathbb{N}$, there is a homomorphism $\phi_i:G \to G_i$ defined by $\phi_{i}(x_k)=x_k$ for $1 \leq k \leq 2i$ and $\phi_{ i}(x_k)=1$ for $k>2i$.

It is straightforward to check that $(G_i, \phi_{i, j})$ is an inverse system with inverse limit $(G, \phi_i)$. (Indeed, let $\phi:G \to \varprojlim G_i$ be the homomorphism induced by the compatible
 homomorphisms $\phi_i:G \to G_i$. The surjectivity of $\phi$ follows from the surjectivity of $\phi_i$ for all $i \in \mathbb{N}$. In order to prove that $\phi$ is injective, it suffices to observe that any homomorphism from $G$ onto a finite $p$-group 
factors through some $\phi_i$.)

Let $H$ be a finitely generated subgroup of $G$. Then there exists $i_0 \in \mathbb{N}$ such that $d(H)=d(\phi_i(H))$ for every $i \geq i_0$. Fix $i \in \mathbb{N}$ with $i >\max\{d(H), i_0\}$.
Then $\phi_i(H)$ has infinite index in $G_i$, and hence it is a free pro-$p$ group (cf. \cite[Chapter I, Section~4.5]{Se97}). 
It follows that ${\phi_i}_{\mid H}: H \to \phi_i(H)$ is an isomorphism, which in turn implies that $H^{ab}$ is torsion-free. 
Therefore, $G$ is locally absolutely torsion-free.

\end{proof}

\section{The subgroups of a direct product}

Goursat's lemma affirms the existence of a one-to-one correspondence between the subgroups of a direct product $G_1 \times G_2$ (of discrete groups) and $5$-tuples $(H_1, N_1, H_2, N_2, \phi)$ where $N_1 \unlhd H_1 \leq G_1$, $N_2 \unlhd H_2 \leq G_2$ and 
$\phi:H_1/N_1 \to H_2/N_2$ is an isomorphism. More precisely, the $5$-tuple $(H_1, N_1, H_2, N_2, \phi)$ determines the subgroup
$K=\{(h_1, h_2) \in H_1 \times H_2 \mid \phi(h_1N_1)=h_2N_2\}$.

In further detail, $K$ is a subdirect product of $H_1$ and $H_2$, i.e., 
$\pi_1(K)=H_1$ and $\pi_2(K)=H_2$, where $\pi_i:G_1 \times G_2 \to G_i$ ($i=1,2$) denote the projection homomorphisms. Moreover,
considering $G_1$ and $G_2$ as subgroups of $G_1 \times G_2$ (in the canonical way), we have $N_1=H_1 \cap K$, $N_2=H_2 \cap K$,  and there is a short exact sequence
\begin{equation}
\label{shortexact1}
\begin{tikzcd}        
   1 \arrow[r] & N_2 \arrow[ r] & K \arrow[r, "\pi_1|_K"] & H_1 \arrow[r] & 1.  
\end{tikzcd}
\end{equation}
By composing $\pi_1|_K$ with the natural epimorphism $H_1 \to H_1/N_1$, we obtain the exact sequence    
\begin{equation}
\begin{tikzcd}        
   1 \arrow[r] & N_1 \times N_2 \arrow[r] & K \arrow[r] & H_1/N_1 \arrow[r] & 1.
\end{tikzcd}
\end{equation}
Therefore, $K/(N_1 \times N_2) \cong H_1/N_1 \cong H_2/N_2$. 

\medskip

Goursat's lemma holds true for pro-$p$ groups -- the proof translates almost verbatim. It goes without saying that when dealing with pro-$p$ groups, subgroups are assumed to be closed and isomorphisms are required to be continuous. To see that $K$ is indeed a closed subgroup of $G_1 \times G_2$, observe that
\[K=\{(h_1, h_2) \in H_1 \times H_2 \mid \psi_1((h_1, h_2))=\psi_2((h_1,h_2))\},\]
where ${\psi_1:H_1 \times H_2 \to H_2/N_2}$  is the composition of the projection $H_1 \times H_2 \to H_1$, the natural epimorphism ${H_1 \to H_1/N_1}$, and the isomorphism $\phi:H_1/N_1 \to H_2/N_2$; and $\psi_2:H_1 \times H_2 \to H_2/N_2$ is the the composition of the projection $H_1 \times H_2 \to H_2$ and the natural epimorphism $H_2 \to H_2/N_2$.


\subsection{Finitely generated subgroups}

There does not seem to exist any straightforward sufficient and necessary condition for a subgroup of a direct product to be finitely generated.  
Be that as it may, the results collected in the following lemma suffice for our purposes.

\begin{lem}
\label{finitely generated}
Let $G_1$ and $G_2$ be pro-$p$ groups, and let $K$ be a subgroup of $G_1 \times G_2$ corresponding to a $5$-tuple $(H_1, N_1, H_2, N_2, \phi)$.
The following assertions hold:
\begin{enumerate}[(i)]
\item If $K$ is finitely generated, then so are $H_1$ and $H_2$.
\item If $H_1$ is finitely generated and $N_2$ is normally finitely generated as a subgroup of $H_2$, then
$K$ is finitely generated.  
\item If $H_1$ and  $H_2$ are finitely generated and $H_1/N_1 \cong H_2/N_2$ is finite, then
$K$ is finitely generated.
\end{enumerate}
\end{lem}

\begin{proof}
$(i)$ follows from the fact that $K$ projects onto $H_1$ and $H_2$. 

\medskip

The proof of $(ii)$ relies on the following short exact sequence, derived from (\ref{shortexact1}):
\begin{equation}
\label{shortexact2}
\begin{tikzcd}        
   1 \arrow[r] & N_2/[H_2, N_2] \arrow[ r] & K/(\{1\} \times [H_2, N_2]) \arrow[r] & H_1 \arrow[r] & 1
\end{tikzcd}
\end{equation}
To obtain a generating set for $K$, first lift a generating set $X=\{x_i \mid i \in I\}$ for $H_1$ to a set $\overline{X}=\{(x_i, x_i') \mid i \in I\} \subseteq K$. 
Next, choose a set $Y=\{y_j \mid j \in J\}$ that normally generates $N_2$ as a subgroup of $H_2$, and set ${\overline{Y}=\{(1, y_j) \mid j \in J\} \subseteq K}$.  
Since $ N_2/[H_2, N_2]$ is generated by the image of $Y$, it follows from (\ref{shortexact2}) that $K/(\{1\} \times [H_2, N_2])$ is generated by the image of 
$Z:=\overline{X} \cup \overline{Y}$. In order to deduce from here that $Z$ generates $K$, it suffices to show that $\{1\} \times [H_2, N_2] \subseteq \Phi(K)$.

For every $h_2 \in H_2$, there is $h_1 \in H_1$ such that $(h_1, h_2) \in K$. Furthermore, given $n_2 \in N_2$, we have $(1, n_2) \in K$ and  
$(1, [h_2, n_2]) = [(h_1, h_2), (1, n_2)] \in \Phi(K)$.
Consequently, $\{1\} \times [H_2, N_2] \subseteq \Phi(K)$ and $Z$ generates $K$. 
Now suppose that $H_1$ is finitely generated and that $N_2$ is normally finitely generated in $H_2$.
Then we may choose $X$ and $Y$ (and thus $Z$) to be finite, which yields $(ii)$.

\medskip

If $H_2$ is finitely generated and  $H_2/N_2$ is finite, then $N_2$ is also finitely generated. Hence, $(iii)$ follows from $(ii)$.
\end{proof}

\begin{rem}
In general, the converse of Lemma~\ref{finitely generated} $(i)$ does not hold. For instance, let $G$ be any finitely generated pro-$p$ group that contains a normal subgroup $N$ that is not normally finitely generated. 
Consider the subgroup $K$ of $G \times G$ corresponding to the $5$-tuple $(G, N, G, N, id_{G/N})$. We claim  that $K$ is not finitely generated. Suppose to the contrary that $X=\{(g_1, g_1'),  \ldots, (g_k, g_k')\} \subseteq K$
is a finite generating set for $K$. For each $1 \leq i \leq k$, there is $n_i \in N$ such that $g_i'=g_in_i$. Let $M$ be the normal closure of $n_1, \ldots, n_k$ in $G$. By assumption, $M$ is a proper subgroup of $N$. Fix $z \in N \setminus M$;
then there is a pro-$p$ word $w(x_1, \ldots, x_k)$ such that 
\[w((g_1, g_1'),  \ldots, (g_k, g_k'))=(w(g_1, \ldots, g_k), w(g_1', \ldots, g_k'))=(1, z).\]
However,
\[z=w(g_1', \ldots, g_k')=w(g_1n_1, \ldots, w_kn_k) \in w(g_1, \ldots, g_k)M=M,\]
which contradicts the choice of $z$.
\end{rem}

\begin{rem}
Although the converse of  Lemma~\ref{finitely generated} $(ii)$ holds under the stronger assumption that $H_1$ is finitely presented (we let the reader provide a proof of this), it is not true in general. 
For an example, let $G$ be a finitely generated {pro-$p$} group that contains a normal subgroup $N$ that is not normally finitely generated.
Set $H:=G \times G/N$, $M:=N \times \{1\}$ and let $\phi:H/M \to H/M$ be the isomorphism defined by 
$\phi((x, yN)M)=(y, xN)M$.   

Consider the subgroup $K$ of $H \times H$ associated to the $5$-tuple $(H, M, H, M, \phi)$. We claim that $K$ is finitely generated.
Since $K / (M \times M) \cong H/M \cong G/N \times G/N$ is a finitely generated pro-$p$ group, it suffices to find a finite set $X \subseteq K$ such that $M \times M \subseteq \overline{\langle X \rangle}$.

Choose a finite generating set $\{g_1, g_2, \ldots, g_k\}$ for $G$. For each $1 \leq i \leq k$, set ${a_i:=((g_i, 1), (1, g_iN))}$, $b_i:=((1, g_iN), (g_i, 1))$,  and let  
$X:=\{a_1, \ldots, a_k, b_1, \ldots, b_k\}$. Given $n \in N$, let $w(x_1, \ldots, x_k)$ be a pro-$p$ word such that $w(g_1, \ldots, g_k)=n$. We have that
\[w(a_1, \ldots, a_k)=((w(g_1, \ldots, g_k), 1), (1, w(g_1, \ldots, g_k )N))=((n, 1), (1, 1))\]
and
\[w(b_1, \ldots, b_k)=((1, w(g_1, \ldots, g_k )N), (w(g_1, \ldots, g_k), 1))=((1, 1), (n, 1)).\]
It follows that $M \times M \subseteq \overline{\langle X \rangle} \subseteq K$, as wanted. 
\end{rem}

\subsection{When $H_2$ is abelian}

\begin{lem}
\label{abelian factor}
Let $G_1$ and $G_2$ be pro-$p$ groups, and let $K$ be a subgroup of $G_1 \times G_2$ corresponding to a $5$-tuple $(H_1, N_1, H_2, N_2, \phi)$.
Suppose that $H_2$ is abelian. Then $[K,K]=[H_1, H_1] \times \{1\}$ and there is a short exact sequence
\[\xymatrix{1 \ar[r] & N_2 \ar[r]^(0.45){\alpha} & K^{ab} \ar[r]^{\beta} & H_1^{ab} \ar[r] & 1  }\] 
where $\alpha(n)=(1, n)[K, K]$ and $\beta((h_1, h_2)[K, K])=h_1[H_1, H_1] $.   
\end{lem}
\begin{proof}
Since $K \leq H_1 \times H_2$,  we have 
$$[K, K] \leq [H_1, H_1] \times [H_2, H_2]=[H_1, H_1] \times \{1\}.$$
For the opposite inclusion, let $h_1, h_1' \in H_1$, and  choose $h_2, h_2' \in H_2$ such that
$\phi(h_1N_1)=h_2N_2$ and $\phi(h_1'N_1)=h_2'N_2$. Then $(h_1, h_2),(h_1', h_2') \in K$, and  since $H_2$ is abelian, we get  
\[([h_1, h_1'], 1)=([h_1, h_1'], [h_2, h_2'])=[(h_1,h_2), (h_1', h_2')] \in [K, K]. \]
It follows that $[H_1, H_1] \times \{1\} \leq [K, K]$, as claimed.

\medskip

It is now clear that $\alpha$ is a monomorphism. Moreover,  $\beta$ is obviously an epimorphism and $\mathrm{Im}(\alpha) \leq \ker(\beta)$; thus 
it only remains to prove that  $\ker(\beta) \leq \mathrm{Im}(\alpha)$.

Let $(h_1, h_2)[K, K] \in\ker(\beta)$. Then $h_1 \in [H_1, H_1]$ and $h_1N_1=N_1$ since $H_1/N_1 \cong H_2/N_2$ is abelian. 
Hence, $h_2N_2=\phi(h_1N_1)=\phi(N_1)=N_2$, and thus, $h_2 \in N_2$. It follows that 
\[(h_1, h_2)[K, K]=(1, h_2)[K, K] \in \mathrm{Im}(\alpha).\] 
\end{proof}

\begin{rem}
\label{split}
If $H_1^{ab}$ is torsion-free (and hence free abelian), the short exact sequence in Lemma~\ref{abelian factor} splits. Hence, in that case, $K^{ab} \cong N_2 \times H_1^{ab}$.
In general, the exact sequence does not split. Indeed, take $G_1=H_1=\mathbb{Z}_p/p\mathbb{Z}_p$, $N_1=\{1\}$, $G_2=H_2=\mathbb{Z}_p$, $N_2=p\mathbb{Z}_p$, and let 
$\phi:H_1 /N_1\to H_2/N_2$ be any isomorphism. Then $N_2 \times H_1^{ab} \cong \mathbb{Z}_p \times \mathbb{Z}_p/p\mathbb{Z}_p$, however, it is easy to see that 
$K^{ab}=K \cong \mathbb{Z}_p$.
\end{rem}

\section{Proof of the Main Theorem}
\subsection{A necessary condition: torsion-free abelianization}

In \cite{SnoTan}, we provided a complete classification of the Frattini-resistant $p$-adic analytic pro-$p$ groups.

\begin{thm}[{cf. \cite[Theorem~1.2]{SnoTan}}]
\label{classification p-adic}
Let $G$ be a $p$-adic analytic pro-p group of dimension $d \geq 1$. Then
$G$ is Frattini-resistant if and only if it is isomorphic to one of the following groups:
\begin{enumerate}[(i)]
\item the abelian group $\mathbb{Z}_p^{d}$;
\item the metabelian group $\langle x \rangle \ltimes \mathbb{Z}_p^{d-1}$, where $\langle x \rangle \cong \mathbb{Z}_p$
and $x$ acts on $\mathbb{Z}_p^{d-1}$ as scalar multiplication by $\lambda $, with $\lambda =1+p^s$ for some $s \geq 1$ if $p > 2$, and $\lambda=1 + 2^s$ for some $s \geq 2$ if $p=2$. 
\end{enumerate}
\end{thm}

The classification of Frattini-resistant $p$-adic analytic pro-$p$ groups will be used in the proof of the second statement of the following lemma. 

\begin{lem}
\label{free abelianization}
Let $G_1$ and $G_2$ be non-trivial pro-$p$ groups. Suppose that $G_1 \times G_2$ is Frattini-resistant.
Then the following assertions hold:
\begin{enumerate}[(i)]
\item Every finitely generated subgroup of $G_1$ or $G_2$ has torsion-free abelianization. 
\item Every $p$-adic analytic subgroup of $G_1$ or $G_2$ is abelian.   
\end{enumerate}  
\end{lem}
\begin{proof}
$(i)$ Let $H_1$ be a finitely generated subgroup of $G_1$. Suppose that $H_1^{ab}$ is not torsion-free. Then there is an element $x \in H_1$ such that
$\bar{x}:=x[H_1, H_1] \in H_1^{ab}$ is a nontrivial element of finite order, say $p^n$, and $H_1^{ab}=\langle \bar{x} \rangle \times A$ for some $A \leq H_1^{ab}$.
Let $N_1:=\pi^{-1}(A)$ where $\pi:H_1 \to H_1^{ab}$ is the natural epimorphism. Then $H_1/N_1$ is a cyclic group of order $p^n$, generated by the image of $x$. 

Being a subgroup of a Frattini-resistant pro-$p$ group (namely, $G_1 \times G_2$), $G_2$ is also Frattini-resistant. 
In particular, $G_2$ is  torsion-free. Choose a nontrivial element $z \in G_2$, and set $H_2:=\langle z \rangle \cong \mathbb{Z}_p$ and $N_2:=H_2^{p^n}=\langle z^{p^n} \rangle$.
Then $H_2/N_2$ is a cyclic group of order $p^n$ and there is an isomorphism $\phi: H_1/N_1 \to H_2/N_2$. 

Consider the subgroup $K$ of $G_1 \times G_2$ associated to the $5$-tuple $(H_1, N_1, H_2, N_2, \phi)$.
By Lemma~\ref{finitely generated} $(iii)$, $K$ is finitely generated.  Let $w:=(x^{p^{n-1}}, 1) \in H_1 \times H_2$.
We have that $w \notin K$, because $\phi(x^{p^{n-1}}N_1) \neq N_2$. On the other hand, it follows from Lemma~\ref{abelian factor} that $w^p=(x^{p^n}, 1) \in [H_1,H_1] \times \{1\}=[K, K] \leq \Phi(K)$, which contradicts the 
Frattini-resistence of $G_1 \times G_2$. 

\medskip
 
Now $(ii)$ follows from $(i)$ and Theorem~\ref{classification p-adic} (the abelianization of 
$\langle x \rangle \ltimes \mathbb{Z}_p^{d-1}$, where $x$ acts as scalar multiplication by $1+p^s$, is 
$\mathbb{Z}_p \times (\mathbb{Z}_p / p^s\mathbb{Z}_p)^{d-1}$).   

\end{proof}

\begin{rem}
It follows from Lemma~\ref{free abelianization} and \cite[Theorem~1.3]{SnoTan} 
that if $G_1 \times G_2$ is Frattini-resistant, then every solvable subgroup of $G_1$ or $G_2$ is abelian.  
\end{rem}

\begin{rem}
Observe that a straightforward adaptation of the proof of Lemma~\ref{free abelianization} shows that
if $G_1 \times G_2$ is strongly Frattini-resistant, then $G_1$ and $G_2$ are absolutely torsion-free.
\end{rem}

\subsection{Finding isomorphic quotients}
A pro-$p$ group $G$ is said to be \emph{powerful} if $p$ is odd and $[G, G] \leq G^p$, or $p=2$ and 
$[G, G] \leq G^{4}$. The reader is referred to \cite{APG} for more details on powerful pro-$p$ groups and 
their role in the theory of {$p$-adic} analytic groups. Apart from the definition, we only need the fact that every finitely generated powerful pro-$p$ group is $p$-adic analytic.

\begin{lem}
\label{isomorphic quotients odd}
Let $p$ be an odd prime, and let $G$ be a pro-$p$ group with $d(G)=2$. If $G$ is not powerful, then  
\[G/(G^p\gamma_3(G)) \cong \langle x_1, x_2 \mid x_1^p, x_2^p, [x_1, x_2]^p,  [[x_1,x_2],x_1], [[x_1,x_2], x_2] \rangle. \]
\end{lem}
\begin{proof}
For notational convenience, set $H:=G/(G^p\gamma_3(G))$ and let $K$ be the pro-$p$ group defined by the above presentation.
Let $\{y_1, y_2\}$ be a generating set for $G$.

Suppose that $[y_1, y_2] \in G^p\gamma_3(G)$. It follows that $[G, G] \leq  G^p\gamma_3(G)$, and thus 
$\Phi(G)=G^p[G, G]=G^p\gamma_3(G)$. In fact, $\Phi(G)=G^p\gamma_k(G)$ for all $k \geq 3$. Indeed, assume by induction that this is true for $k=i \geq 3$; then
\begin{align*}
\gamma_{i}(G)&=[\gamma_{i-1}(G), G] \leq [\Phi(G),G]=[G^p\gamma_i(G), G] \\ 
&=[G^p, G]\gamma_{i+1}(G) \leq G^p\gamma_{i+1}(G).  
\end{align*}
Hence, $\Phi(G)=G^p\gamma_{i}(G)=G^p\gamma_{i+1}(G)$.

It follows that 
\[\Phi(G)=\bigcap_{i=1}^{\infty} G^p\gamma_i(G)=G^p \bigg(\bigcap_{i=1}^{\infty}\gamma_{i}(G) \bigg)=G^p,
\]
which implies that $G$ is powerful.

Now assume that $G$ is not powerful. Then as we just observed $[y_1, y_2] \notin G^p\gamma_3(G)$.
Moreover, it is easy to see that $[H, H]$ is a cyclic group of order $p$, generated by the image of $[y_1, y_2]$. 
It is also clear that $[K, K]=\langle [x_1, x_2] \rangle$ is a cyclic group of order at most $p$.
Furthermore, $H^{ab} \cong G/\Phi(G)$ and $K^{ab}$ are elementary abelian $p$-groups of order $p^2$.
Hence, $|K| \leq |H|=p^3$. 

Now consider the epimorphism $\phi: K \to H$ defined by $\phi(x_1)=y_1$ and $\phi(x_2)=y_2$.
Since $|K|\leq |H|$, $\phi$ must be an isomorphism.
\end{proof}

For $p=2$, there is a similar result.

\begin{lem}
\label{isomorphic quotients 2}
Let $G$ be a pro-$2$ group with $d(G)=2$. Assume that $G^{ab}$ does not contain a direct cyclic factor of order $2$.
If $G$ is not powerful, then 
\[G/(G^4\gamma_3(G)) \cong \langle x_1, x_2 \mid x_1^4, x_2^4, [x_1, x_2]^2,  [[x_1,x_2],x_1], [[x_1,x_2], x_2] \rangle. \]
\end{lem}
\begin{proof}
Put $H:=G/(G^4\gamma_3(G))$ and let $K$ be the pro-$2$ group defined by the above presentation.
Fix a generating set $\{y_1, y_2\}$ for $G$. 

Suppose that $[y_1, y_2] \in G^4\gamma_3(G)$. Consequently, $[G, G] \leq  G^4\gamma_3(G)$. 
We will show that this implies that $[G, G] \leq G^4\gamma_k(G)$ for every $k \geq 3$. 
Indeed, assume by induction that this is true for $k=i \geq 3$; then
\begin{align*}
\gamma_{i}(G)&=[\gamma_{i-1}(G), G] \leq [[G, G] ,G] \leq [G^4\gamma_i(G), G] \\ 
&=[G^4, G]\gamma_{i+1}(G) \leq G^4\gamma_{i+1}(G).  
\end{align*}
Hence, $[G, G] \leq G^4\gamma_{i}(G)=G^4\gamma_{i+1}(G)$.

It follows that 
\[[G, G] \leq \bigcap_{i=1}^{\infty} G^4\gamma_i(G)=G^4 \bigg(\bigcap_{i=1}^{\infty}\gamma_{i}(G) \bigg)=G^4,
\]
which implies that $G$ is powerful.

Assume that $G$ is not powerful. Then $[y_1, y_2] \notin G^4 \gamma_3(G)$.
Furthermore, since $[y_1, y_1] \in G^2$, we have $[y_1, y_1]^2 \in G^4 \leq G^4\gamma_3(G)$.
Therefore, $[H, H]$ is a cyclic group of order $2$, generated by the image of $[y_1, y_2]$. It follows from the assumption on the abelianization of $G$ that $H^{ab} \cong \mathbb{Z}_2 /4 \mathbb{Z}_2 \times \mathbb{Z}_2 /4\mathbb{Z}_2$.
Hence, $H$ has order $32$.  

On the other hand, $[K, K]=\langle [x_1, x_2] \rangle$ is a cyclic group of order at most $2$,
and $K^{ab} \cong \mathbb{Z}_2 /4 \mathbb{Z}_2 \times \mathbb{Z}_2 /4\mathbb{Z}_2$. 
Therefore, $|K| \leq |H|$. 

Now consider the epimorphism $\phi: K \to H$ defined by $\phi(x_1)=y_1$ and $\phi(x_2)=y_2$.
Since $|K|\leq |H|$, $\phi$ must be an isomorphism.
\end{proof}

\begin{rem}
\label{nice iso}
We record an important consequence of the proofs of Lemma~\ref{isomorphic quotients odd} and Lemma~\ref{isomorphic quotients 2}. Let $G_1$ and $G_2$ be non-powerful pro-$p$ groups with $d(G_1)=d(G_2)=2$. In addition, if $p=2$, assume that $G_1^{ab}$ and $G_2^{ab}$ do not contain direct cyclic factors of order $2$.
Fix generating sets $\{x_1, y_1\}$ and $\{x_2, y_2\}$ for $G_1$ and $G_2$, respectively. Then $[x_1, y_1] \notin G_1^p\gamma_3(G_1)$ and  $[x_2, y_2] \notin G_2^p\gamma_3(G_2)$ if $p$ is odd, and 
$[x_1, y_1] \notin G_1^4\gamma_3(G_1)$ , $[x_2, y_2] \notin G_2^4\gamma_3(G_2)$, $[x_1, y_1]^2 \in G_1^4\gamma_3(G_1)$ and $[x_2, y_2]^2 \in G_2^4\gamma_3(G_2)$ if $p=2$. Moreover,  
there is an ismorphism ${\phi:G_1/G_1^p\gamma_3(G_1) \to G_2/G_2^p\gamma_3(G_2)}$  such that $\phi(\bar{x}_1)=\bar{x}_2$ and $\phi(\bar{y}_1)=\bar{y}_2$ if $p$ is odd, and an isomorphism  
${\phi:G_1/G_1^4\gamma_3(G_1) \to G_2/G_2^4\gamma_3(G_2)}$  such that $\phi(\bar{x}_1)=\bar{x}_2$ and $\phi(\bar{y}_1)=\bar{y}_2$ if $p=2$. 
\end{rem}

\subsection{The final step}
\begin{proof}[Proof of Theorem~\ref{main theorem}]
We prove $(i)$, the proof of $(ii)$ is similar (even simpler).
Suppose that all finitely generated subgroups of $G_1$ have torsion-free abelianization and that $G_2$ is torsion-free abelian.
It follows from Lemma~\ref{finitely generated} $(i)$ and Lemma~\ref{abelian factor} (see also Remark~\ref{split}) that all finitely generated subgroups of ${G_1 \times G_2}$ have torsion-free abelianization.
Therefore, $G_1 \times G_2$ is commutator-resistant, and thus it is Frattini-resistant. 

\medskip

For the converse, suppose that $G_1 \times G_2$ is a Frattini-resistant pro-$p$ group. It follows from Lemma~\ref{free abelianization} $(i)$ that all finitely generated subgroups of $G_1$ and $G_2$ have torsion-free abelianization.
Therefore, it remains to prove that at least one of the groups $G_1$ or $G_2$ is abelian. 

Suppose to the contrary that there are elements $x_1, y_1 \in G_1$ and $x_2, y_2 \in G_2$ such that $[x_1, y_1] \neq 1$ and $[x_2, y_2] \neq 1$.  
Let $H_1:=\langle x_1, y_1 \rangle$ and $H_2:=\langle x_2, y_2 \rangle$. 
Set $N_1:=H_1^p\gamma_3(H_1)$ and $N_2:=H_2^p\gamma_3(H_2)$ if $p$ is odd, and $N_1:=H_1^4\gamma_3(H_1)$ and $N_2:=H_2^4\gamma_3(H_2)$ if $p=2$.
By Lemma~\ref{free abelianization} $(ii)$, neither $H_1$ nor $H_2$ is powerful. It follows from Lemma~\ref{isomorphic quotients odd} and Lemma~\ref{isomorphic quotients 2}, 
that $H_1/N_1$ and $H_2/ N_2$ are isomorphic \mbox{pro-$p$} groups. More precisely, there is an isomorphism $\phi:H_1/ N_1 \to H_2/N_2$ satisfying $\phi(x_1N_1)=x_2N_2$ and 
$\phi(y_1N_1)=y_2N_2$ (see Remark~\ref{nice iso}).

Let $K$ be the subgroup of $G_1 \times G_2$ determined by the $5$-tuple $(H_1, N_1, H_2, N_2, \phi)$. Since $H_1/N_1 \cong H_2/ N_2$  is a finite group,
it follows from Lemma~\ref{finitely generated} $(iii)$ that $K$ is finitely generated.

Clearly, $\gamma_3(H_1) \times \{1\} \leq K$. Consequently, 
\[\gamma_4(H_1) \times \{1\}=[H_1, \gamma_3(H_1)] \times \{1\}=[K, \gamma_3(H_1) \times \{1\}] \leq [K, K].\]
Now the proofs for $p$ an odd prime and $p=2$ diverge.  

\medskip

\textbf{Case 1:} \textit{$p$ is an odd prime.} By a straightforward calculation,  
$$[x_1, [x_1, y_1]]^p \gamma_4(H_1)=[{x_1}^p, [x_1, y_1]] \gamma_4(H_1).$$ 
Since $([x_1, [x_1, y_1]], 1)$, $([x_1, y_1], [x_2, y_2])$ and $(x_1^p, 1)$ all belong to $K$, we get 
\begin{align*}
([x_1, [x_1, y_1]], 1)^p[K, K]&=([x_1, [x_1, y_1]]^p, 1)[K, K]=([x_1^p, [x_1, y_1]], 1)[K, K] \\
&=[({x_1}^p, 1), ([x_1, y_1], [x_2, y_2])][K, K]=[K, K].
\end{align*}

Consider the element $w:=([x_1, y_1], 1) \in H_1 \times H_2$. By Remark~\ref{nice iso}, $[x_1, y_1] \notin N_1$, and thus $w \notin K$.
On the other hand, $w^p=([x_1, y_1]^p, 1) \in K$. In fact, we contend that $w^p \in [K, K]$, which would yield a contradiction with the assumption that $G_1 \times G_2$ is Frattini-resistant. 

To see this, first observe that 
\begin{align*}
[x_1^p, y_1]\gamma_4(H_1)&=[x_1, y_1]^{{x_1}^{p-1}}[x_1, y_1]^{{x_1}^{p-2}} \cdots [x_1, y_1]\gamma_4(H_1)\\
&=[x_1, y_1]^p[[x_1, y_1], x_1]^{p(p-1)/2}\gamma_4(H_1).
\end{align*}

Using that $p$ is an odd prime, $\gamma_4(H_1) \times \{1\} \leq [K, K]$ and ${([x_1, [x_1, y_1]], 1)^p \in [K, K]}$, we get
\begin{align*}
w^p[K, K]&=([x_1, y_1]^p, 1)[K, K]=([x_1^p, y_1], 1)([x_1,[x_1, y_1], 1)^{p(p-1)/2}[K, K]\\
&=([x_1^p, y_1], 1)[K, K]=[(x_1^p, 1), (y_1, y_2)][K, K]=[K, K].
\end{align*}
Therefore, $w^p \in [K, K]$, as claimed. 

\medskip

\textbf{Case 2:} \textit{$p=2$.} Put $w:=([x_1, y_1], 1) \in H_1 \times H_2$. By Remark~\ref{nice iso}, $w \notin K$ but $w^2 \in K$.
Since $G_1 \times G_2$ is Frattini-resistant, it follows that $w^2 \notin \Phi(K)$. We will show that that $w^4 \in [K, K]$, which in light of Theorem~\ref{commutator-resistant}, 
would yield a contradiction. Indeed, 
\begin{align*}
[x_1^4, y_1] \gamma_4(H_1)&=[x_1, y_1]^{{x_1}^{3}}[x_1, y_1]^{{x_1}^{2}}[x_1, y_1]^{{x_1}}[x_1, y_1]\gamma_4(H_1)\\
&=[x_1, y_1]^4[[x_1, y_1], x_1]^6\gamma_4(H_1)=[x_1, y_1]^4[[x_1, y_1]^6, x_1]\gamma_4(H_1),
\end{align*}
which implies that
\begin{align*}
w^4[K, K]&=([x_1, y_1]^4, 1)[K, K]=([x_1^4, y_1], 1)([x_1, [x_1, y_1]^6], 1)[K, K]\\
&=[(x_1^4, 1), (y_1, y_1)] \cdot [(x_1, x_2), ([x_1, y_1]^6, 1)][K, K]=[K, K].
\end{align*}
\end{proof}

\section{Proofs of the Corollaries}

\begin{proof}[Proof of Corollary~\ref{product of Galois groups}]

First assume that $\sqrt{-1} \in k$ if $p=2$. By \cite[Theorem~7.1 and Theorem~7.2]{SnoTan}, $G_k(p)$ is strongly Frattini-resistant, and  it follows from Theorem~\ref{main theorem} that one of the groups $G_1$ or $G_2$ is abelian and the other one is absolutely torsion-free. 

Put $F:=k(\mu_{\infty})$, $H:=\mathrm{Gal}(\sqrt[\infty]{k}/k)$ and $A:=\mathrm{Gal}(\sqrt[\infty]{k}/F)$, where $\sqrt[\infty]{k}$ denotes the extension of $k$ generated by all roots of order a power of $p$ of elements in $k$ (taken inside $k(p)$). It is readily seen that the restriction homomorphism ${G_p(k) \twoheadrightarrow H}$ is a Frattini cover (i.e. its kernel is contained in the Frattini subgroup of $G_k(p)$); hence  $d(H)=d(G_k(p)) \geq 2$. 
Since  $H/A \cong \mathrm{Gal}(F/k)$ is a cyclic pro-$p$ group, it follows that $A \neq 1$. Moreover, by Kummer theory, $A$  is torsion-free abelian (cf. \cite[Theorem~4.2]{EQ1}). 

Let $\theta: G_k(p) \to 1+p\mathbb{Z}_p$ be the cyclotomic pro-$p$ character defined by 
$\sigma(\zeta)=\zeta^{\theta (\sigma)}$ for all $\sigma \in G_k(p)$ and $\zeta \in \mu_{\infty}$. 
Given $\sigma \in H$, $\tau \in A$ and $\sqrt[p^s]{a} \in \sqrt[\infty]{k}$ (for some $a \in k$ and $s \geq 1$), 
write $\sigma^{-1}(\sqrt[p^s]{a})=\zeta^m\sqrt[p^s]{a}$ and $\tau(\sqrt[p^s]{a})=\zeta^n\sqrt[p^s]{a}$, where $\zeta$ denotes a fixed primitive $p^s$th root of unity. 
We have 
\begin{equation}
\label{conjugation}
\begin{split}
\sigma\tau\sigma^{-1}(\sqrt[p^s]{a})&=\sigma\tau(\zeta^m\sqrt[p^s]{a})=\sigma(\zeta^{m+n}\sqrt[p^s]{a})=\sigma(\zeta^{n})\sigma(\zeta^m\sqrt[p^s]{a})\\
&=\zeta^{n\theta (\sigma)}\sqrt[p^s]{a}=\tau^{\theta(\sigma)}(\sqrt[p^s]{a}).
\end{split}
\end{equation}
Hence, $\sigma\tau\sigma^{-1}=\tau^{\theta(\sigma)}$. 

Since one of the direct factors of $G_k(p)$ is abelian, there is a central element $\rho$ in $G_k(p)$ that does not belong to $\Phi(G_k(p))$.
The Frattini cover $G_p(k) \twoheadrightarrow H$ maps $\rho$ to a non-trivial central element in $H$. The calculation (\ref{conjugation}) now implies that the cyclotomic character $\theta$ must be trivial. 
(To see this, first let $\rho$ take the place of $\sigma$ and conclude that $\theta(\rho)=1$, or equivalently, $\rho \in A$;
next let $\rho$ take on the role of $\tau$.) Therefore, $\mu_{\infty} \subseteq k$. 

It remains to demonstrate the impossibility of $\sqrt{-1} \notin k$ when $p=2$. By \cite[Lemma~2.1]{Koen}, $G_k(2)$ is torsion-free.   
Without loss of generality, we may assume that there is an automorphism $\sigma \in G_1$ such that $\sigma(\sqrt{-1})=-\sqrt{-1}$. Choose a non-trivial element $\tau \in G_2$, and
consider the subgroup $T:= \langle \sigma \rangle \times \langle \tau \rangle \cong \mathbb{Z}_2 \times \mathbb{Z}_2$ of $G_k(2)$. Clearly, $T$ is a maximal pro-$p$ Galois group of some field that does not contain $\sqrt{-1}$.
Hence, we may assume from the outset that $G_1$ and $G_2$ are infinite cyclic pro-$2$ group.

Let $M=\mathrm{Gal}(k(p)/k(\sqrt{-1}))$,  $N_1:=M \cap G_1$, $N_2:=M \cap G_2$, and let $L$ be the fixed field of $N_1 \times N_2$.
It follows from what we have already proved that $\mu_{\infty} \subseteq L$. In fact, $\mu_{\infty} \subseteq k(\sqrt{-1})$, 
since otherwise $L$ would be an extension of $k(\sqrt{-1})$ of degree 2 generated by some primitive $2^s$th root of unity $\zeta$ (for some $s \geq 8$)  with $\zeta^2 \in k(\sqrt{-1})$.
However, taking norms into $k(\sqrt{-1})$, we get $N(\sqrt{\zeta})^2=N(\zeta)=-\zeta^2$, which yields a contradiction with $\zeta \notin k(\sqrt{-1})$.

Hence, we have that $\mathrm{Gal}(k(\mu_{\infty})/k)=\mathrm{Gal}(k(\sqrt{-1})/k)$ is a cyclic group of order $2$, and 
we can argue as before that $A \neq 1$ and that the cyclotomic character must be trivial. Therefore, $\mu_{\infty} \subseteq k$, which contradicts our assumption.  
 
\end{proof}

\begin{proof}[Proof of Corollary~\ref{smooth}]
If one of the groups $G_1$ or $G_2$ is abelian and the other one is absolutely torsion-free, then $G_1 \times G_2$ is absolutely torsion free.  
Consequently, $G_1 \times G_2$ is $1$-smooth with respect to the trivial homomorphism from $G_1 \times G_2$ to $1+p\mathbb{Z}_p$.

Conversely, suppose that $(G_1 \times G_2, \theta)$ is $1$-smooth (with $\textrm{Im}(\theta) \leq 1 + 4\mathbb{Z}_2$ if $p=2$).
By \cite[Theorem~1.11]{SnoTan}, $G_1 \times G_2$ is strongly Farttini-resistant, and  it follows from Theorem~\ref{main theorem} that one of the groups $G_1$ or $G_2$ is abelian and the other one is 
absolutely torsion-free. Furthermore, \cite[Theorem~1.11]{SnoTan} implies that $\theta$ must be the trivial homomorphism.
\end{proof}
    
\begin{proof}[Proof of Corollary~\ref{not strongly Frattini-resistant}]
By Proposition~\ref{Demushkin of inf rank}, $G$ is locally absolutely torsion free but not absolutely torsion free.
Hence, it follows from Theorem~\ref{main theorem} that $G \times \mathbb{Z}_p$ is Frattini-resistant but not strongly Frattini-resistant.
\end{proof}


\begin{thebibliography}{10}
\bibitem{CleFlo17} C.~De Clercq and M.~Florence,  \textit{Lifting theorems and smooth profinite groups}, preprint, available
at arXiv:1710.10631, 2017.
\bibitem{APG} J.~D.~Dixon, M.~P.~F.~du Sautoy, A.~Mann, and D.~Segal, Analytic pro-p groups, Cambridge
Studies in Advanced Mathematics 61, Cambridge University Press, Cambridge, second edition, (1999).
\bibitem{EQ1} I.~Efrat, C.~Quadrelli, \textit{The Kummerian property and maximal pro-$p$ Galois groups}, J. Algebra, \textbf{525} (2019), 284 -- 310. 
\bibitem{Koen} J.~Koenigsmann, \textit{Products of absolute Galois groups}, International Mathematics Research Notices, \textbf{2005} (2005), no. 24, 1465 -- 1486.
\bibitem{QuWe18} C.~Quadrelli and Th.~Weigel,  \textit{Profinite groups with a cyclotomic $p$-orientation}, Documenta Math.  \textbf{25} (2020), 1881--1916. 
\bibitem{Ribes1} L.~Ribes and P.~Zalesskii, Profinite groups. Second edition. Springer-Verlag, Berlin (2010).
\bibitem{Se97} J.~P.~Serre, Galois Cohomology. Springer-Verlag, Berlin (1997). 
\bibitem{SnoTan} I.~Snopce and S.~Tanushevski, \textit{Frattini-injectivity and maximal pro-$p$ Galois groups}, arXiv:2009.09297.
\bibitem{Wurfel} T.~W{\"u}rfel, \textit{On a class of pro-$p$ groups occurring in Galois theory},  J. Pure Appl. Algebra \textbf{36}
(1985), 95 -- 103.
 
\end{thebibliography}
\end{document}